\documentclass[12pt]{amsart}
\usepackage{latexsym}
\usepackage{amscd,amssymb,times,fullpage}
\newtheorem{thm}{Theorem}[section]
\newtheorem{prop}{Proposition}[section]
\newtheorem{lem}{Lemma}[section]
\newtheorem{cor}{Corollary}[section]
\newtheorem{conj}{Conjecture}[section]

\theoremstyle{remark}
\newtheorem{rem}{Remark}[section]

\theoremstyle{definition}

\title{On the conformal systoles of four-manifolds}
\author{M.~J.~D.~Hamilton}
\address{Mathematisches Institut, Ludwig-Maximilians-Universit\"at M\"unchen,
Theresienstr.~39, 80333 M\"unchen, Germany}
\email{Mark.Hamilton@mathematik.uni-muenchen.de}
\date{November 16, 2006; MSC 2000: primary 53C23, 53D35, secondary 57M50, 57R57}

\begin{document}

\begin{abstract} We extend a result of M.~Katz on conformal systoles to all four-manifolds with $b_2^+=1$ which have odd intersection form. The same result holds for all four-manifolds with $b_2^+=1$ with even intersection form and which are symplectic or satisfy the so-called $\frac{5}{4}$-conjecture. 
\end{abstract}

\maketitle

\section{Introduction}

There are several notions of systolic invariants for Riemannian manifolds, which were introduced by M.~Berger and M.~Gromov (see \cite{G} and \cite{B,CK} for an overview). The most basic concept is the $k$-systole $\mbox{sys}_k(X,g)$ of a Riemannian manifold $X$, defined as the infimum over the volumes of all cycles representing non-zero classes in $H_k(X;\mathbb{Z})$. In this note we discuss a different systole, namely the {\it conformal systole}, which depends only on the conformal class of the Riemannian metric. We briefly review its definition (see Section \ref{def} for details). 

Let $(X^{2n},g)$ be a closed oriented even dimensional Riemannian manifold. The Riemannian metric defines an $L^2$-norm on the space of harmonic $n$-forms on $X$ and hence induces a norm on the middle-dimensional cohomology $H^n(X;\mathbb{R})$. The {\it conformal $n$-systole} $\mbox {confsys}_n(X,g)$ is the smallest norm of a non-zero element in the integer lattice $H^n(X;\mathbb{Z})_{\mathbb{R}}$ in $H^n(X;\mathbb{R})$. It is known that for a fixed manifold $X$ the conformal $n$-systoles are bounded from above as $g$ varies over all Riemannian metrics. Hence the supremum $CS(X)=\sup_g \mbox{confsys}_n(X,g)$ of the conformal systoles over all metrics $g$ is a finite number, which is {\sl a priori} a diffeomorphism invariant of $X$.    

The interest in the literature has been to find bounds for $CS(X)$ that depend only on the topology of $X$, e.g.~the Euler characteristic of $X$, where $X$ runs over some class of manifolds. In \cite{BS} P.~Buser and P.~Sarnak proved the following inequalities for the closed orientable surfaces $\Sigma_s$ of genus $s$: there exists a constant $C>0$ independent of $s$ such that  
\begin{equation}
C^{-1}\log s < CS(\Sigma_s)^2<C\log s,\quad \forall s\geq 2.
\end{equation}

In dimension 4, M.~Katz \cite{Ka} proved a similar inequality for the conformal 2-systole of blow-ups of the complex projective plane $\mathbb{C}P^2$: there exists a constant $C>0$ independent of $n$ such that 
\begin{equation}
C^{-1}\sqrt{n}<CS(\mathbb{C}P^2 \# n\overline{\mathbb{C}P^2})^2 <C n, \quad  \forall n>0.
\end{equation}
In his proof, M.~Katz used a conjecture on the {\it period map} of 4-manifolds $X$ with $b_2^+=1$. The period map is defined as the map taking a Riemannian metric $g$ to the point in the projectivization of the positive cone in $H^2(X;\mathbb{R})$ given by the $g$-selfdual direction (see Section \ref{def}). The conjecture, which is still open, claims that this map is surjective. However, an inspection of the proof of M.~Katz shows that this surjectivity conjecture in full strength is not needed and that in fact his theorem holds in much greater generality.

In Section \ref{proofs}, we first remark that the following proposition holds as a consequence of recent work of D.~T.~Gay and R.~Kirby \cite{GK}.
\begin{prop} \label{dense} The period map for all closed 4-manifolds with $b_2^+=1$ has dense image.
\end{prop}

Using the argument of M.~Katz, this implies the following theorem.

\begin{thm} \label{katz} There exists a universal constant $C$ independent of $X$ and $n=b_2(X)$ such that
\begin{equation}
C^{-1}\sqrt{n} < CS(X)^2 <C n,
\end{equation}
for all closed 4-manifolds $X$ with $b_2^+=1$ which have odd intersection form.
\end{thm}
 
Another consequence of Proposition \ref{dense} is the following theorem.

\begin{thm} \label{coarse} Let $X,X'$ be closed 4-manifolds with $b_2^+=1$ which have isomorphic intersection forms. Then $CS(X)=CS(X')$. 
\end{thm}
This shows that in dimension 4 the invariant $CS$ is much coarser than a diffeomorphism invariant. Theorem \ref{coarse} can be compared to a result of I.~K.~Babenko (\cite{Ba}, Theorem 8.1.), who showed that a certain 1-dimensional systolic invariant for manifolds of arbitrary dimension is a homotopy-invariant. 

Theorem \ref{coarse} enables us to deal with even intersection forms. Suppose $X$ is a closed 4-manifold with $b_2^+=1$ and even intersection form. By the classification of indefinite even quadratic forms, the intersection form of $X$ is isomorphic to $H\oplus (-k)E_8$ for some $k\geq 0$. In particular, for each $r\in \mathbb{N}$ there are only finitely many possible even intersection forms of rank less or equal than $r$. Hence by Theorem \ref{coarse}, the invariant $CS$ takes only finitely many values on all 4-manifolds with even intersection form, $b_2^+=1$ and $b_2\leq r$.  We will show that {\it symplectic} 4-manifolds $X$ with $b_2^+(X)=1$ and even intersection form necessarily have $b_2(X)\leq 10$ (see Section \ref{symplectic}). The same bound holds if $X$ satisfies the so-called {\it $\frac{5}{4}$-conjecture} (see Section \ref{5/4}). Hence together with Theorem \ref{katz}, we get the following corollary, which possibly covers all 4-manifolds with $b_2^+=1$. 

\begin{cor} There exists a  universal constant $C$ independent of $X$ and $n=b_2(X)$ such that 
\begin{equation}
C^{-1}\sqrt{n} < CS(X)^2 < C n, 
\end{equation}
for all closed 4-manifolds $X$ with $b_2^+=1$ which are symplectic or have odd intersection form $Q$ or satisfy the $\frac{5}{4}$-conjecture if $Q$ is even. 
\end{cor}  

\section{Definitions}\label{def}
Let $(X^{2n},g)$ be a closed oriented Riemannian manifold. We denote the space of $g$-harmonic $n$-forms on $X$ by $\mathcal{H}^n(X)$. The Riemannian metric defines an $L^2$-norm on $\mathcal{H}^n(X)$ by
\begin{equation}
|\alpha|^2_{L^2}  =\int_X\alpha\wedge *\alpha,\quad \alpha\in \mathcal{H}^n(X),
\end{equation} 
where $*$ is the Hodge operator.

Given the unique representation of cohomology classes by harmonic forms, we obtain an induced norm $| \cdot |_g$, which we call the 
{\it $g$-norm}, on the middle-dimensional cohomology $H^n(X;\mathbb{R})$. The conformal $n$-systole is defined by 
\begin{equation}
\mbox{confsys}_n(X,g) = \min\{|\alpha|_g\mid \alpha\in H^n(X;\mathbb{Z})_\mathbb{R}\setminus \{0\}\},
\end{equation} 
where $H^n(X;\mathbb{Z})_\mathbb{R}$ denotes the integer lattice in $H^n(X;\mathbb{R})$. More generally, if $L$ is any lattice with a norm $|\cdot |$, we define 
\begin{equation}
\lambda_1(L,|\cdot |)=\min\{|v|\mid v\in L\setminus\{0\}\},
\end{equation}
hence $\mbox{confsys}_n(X,g)=\lambda_1(H^n(X;\mathbb{Z})_\mathbb{R}, |\cdot|_g)$. The conformal systole depends only on the conformal class of $g$ since the Hodge star operator in the middle dimension is invariant under conformal changes of the metric.

The conformal systoles satisfy the following universal bound (see \cite{Ka} equation (4.3)):
\begin{equation}\label{bound}
\mbox{confsys}_n(X,g)^2 < \frac{2}{3}b_n(X), \quad\mbox{for}\,\, b_n(X)\geq 2.
\end{equation}
Clearly, there is also a bound for $b_n(X)=1$, since the Hodge operator on harmonic forms is up to a sign the identity in this case, hence $\mbox{confsys}_n(X,g)=1$. Therefore, the supremum 
\begin{equation}
CS(X)=\sup_g\mbox{confsys}_n(X,g)
\end{equation}
is well-defined for all closed orientable manifolds $X^{2n}$.

We now consider the case of 4-manifolds, $n=2$. In this case the $g$-norm on $H^2(X;\mathbb{R})$ is related to the intersection form $Q$ by the following formula:
\begin{equation}
|\alpha|_g^2=Q(\alpha^+,\alpha^+)-Q(\alpha^-,\alpha^-),
\end{equation}
where $\alpha=\alpha^++\alpha^-$ is the decomposition given by the splitting $H^2(X;\mathbb{R})=H^+\oplus H^-$ into the subspaces represented by $g$-selfdual and anti-selfdual harmonic forms. We abbreviate this formula to
\begin{equation}
|\cdot|_g^2=SR(Q,H^-),
\end{equation}
where $SR$ denotes sign-reversal. Since $H^+$ is the $Q$-orthogonal complement of $H^-$, we conclude that the norm $|\cdot|_g$ is completely determined by the intersection form and the $g$-anti-selfdual subspace $H^-$. 

In particular, let $X$ be a closed oriented 4-manifold with $b_2^+=\dim H^+=1$. The map which takes a Riemannian metric to the selfdual line $H^+$ in the cone $\mathcal{P}$ of elements of positive square in $H^2(X;\mathbb{R})$ (or to the point in the projectivization $\mathbb{P}(\mathcal{P})$ of this cone) is called the period map. In the proof of his theorem, M.~Katz used the following conjecture, which is still open, in the case of blow-ups of $\mathbb{C}P^2$.  
\begin{conj} The period map is surjective for all closed oriented 4-manifolds with $b_2^+=1$.
\end{conj}
If $X$ is a 4-manifold with $b_2^+=1$, we switch the orientation (this does not change the $g$-norm on $H^2(X;\mathbb{R})$) to obtain $b_2^-=1$. Then the $g$-norm is completely determined by the intersection form and the selfdual line, which in the new orientation is $H^-$.
\begin{lem} \label{cont} Let $\bar{X}$ be a 4-manifold with $b_2^-=1$ and intersection form $\bar{Q}$ and let $L$ be the integer lattice in $H^2(\bar{X};\mathbb{R})$. Then $\lambda_1(L,SR(\bar{Q},V)^{1/2})$ depends continuously on the anti-selfdual line $V$.
\end{lem}
This follows because the vector space norm $SR(\bar{Q},V)^{1/2}$ depends continuously on $V$ and the minimum in $\lambda_1$ cannot jump (cf.~Remark 9.1.~in \cite{Ka}). 

\section{Proof of Proposition \ref{dense}, Theorem \ref{katz} and Theorem \ref{coarse}} \label{proofs}

The following theorem can be deduced from recent work of D.~T.~Gay and R.~Kirby \cite{GK} (compare also \cite{ADK}). 
\begin{thm} If $X$ is a closed oriented 4-manifold and $\alpha\in H^2(X;\mathbb{Z})_\mathbb{R}$ a class of positive square, then there exists a Riemannian metric on $X$ such that the harmonic representative of $\alpha$ is selfdual.
\end{thm}
This implies Proposition \ref{dense}, because the set of points given by the lines through integral classes in $H^2(X;\mathbb{R})$ form a dense subset of $\mathbb{P}(\mathcal{P})$. We can now prove Theorem \ref{coarse}.
\begin{proof} Let $\bar{X}$ be $X$ with the opposite orientation, $L$ be the integer lattice in $H^2(\bar{X};\mathbb{R})$ and $\bar{Q}=-Q$. We have 
\begin{equation} CS(X) \leq \sup_V \lambda_1(L,SR(\bar{Q},V)^{1/2}),
\end{equation}
where the supremum extends over all negative definite lines $V$ in $H^2(\bar{X};\mathbb{R})$. This inequality is an equality because the image of the period map is dense and because of Lemma \ref{cont}. The right-hand side depends only on the intersection form.
\end{proof}

We now prove Theorem \ref{katz}. 
\begin{proof} Let $X$ be a closed 4-manifold with intersection form $Q \cong (1)\oplus n(-1)$ for some $n>0$. It is enough to prove inequalities of the form $A\sqrt{n}<CS(X)^2<Bn$ for some constants $A,B>0$, since we can then take $C=\max\{A^{-1},B\}$. The inequality on the right-hand side follows from equation (\ref{bound}). We are going to prove the inequality on the left, following the proof of M.~Katz. 
\begin{lem} \label{conway} There exists a constant $k(n)>0$ (which depends only on $n$ and is asymptotic to $n/2\pi e$ for large $n$) such that 
\begin{equation}
CS(X) \geq k(n)^{1/4}.
\end{equation}
\end{lem}
\begin{proof}
It is more convenient to work with $\bar{X}$, which is $X$ with the opposite orientation. We identify $H^2(\bar{X};\mathbb{R})=\mathbb{R}^{n,1}=\mathbb{R}^{n+1}$ with quadratic form $q_{n,1}$ given by $q_{n,1}(x)=x_1^2+...+x_n^2-x_{n+1}^2$. If $g$ is a metric on $\bar{X}$ then the $g$-norm is given by $|\cdot|^2_g=SR(q_{n,1},v)$ where $*v=-v$ and $SR$ means sign reversal in the direction of $v$. 

Let $L=I_{n,1}=\mathbb{Z}^{n,1}\subset \mathbb{R}^{n,1}$ be the integer lattice. According to Conway-Thompson (see \cite{MH}, Ch.~II, Theorem 9.5), there exists a positive definite odd integer lattice $CT_n$ of rank $n$ with 
\begin{equation}
\min_{x\in CT_n \setminus\{0\}} x\cdot x \geq k(n),
\end{equation}
where $k(n)$ is asymptotic to $n/2\pi e$ for $n\rightarrow \infty$. By the classification of odd indefinite unimodular forms, $CT_n\oplus I_{0,1}\cong I_{n,1}$, hence there exists a vector $v \in \mathbb{Z}^{n,1}$ with $q_{n,1}(v)=-1$ such that $v^\perp \cong CT_n$.

According to M.~Katz, there exists an isometry $A$ of $(\mathbb{R}^{n,1},q_{n,1})$ such that \begin{equation}
\lambda_1(L,SR(q_{n,1},Av)^{1/2})\geq k(n)^{1/4}.
\end{equation}
By Proposition \ref{dense}, there exists a sequence of Riemannian metrics $g_i$ on $X$ whose selfdual lines converge to the line through $Av$. Lemma \ref{cont} implies 
\begin{equation}
\mbox{confsys}_n(X,g_i)\stackrel{i\rightarrow \infty}{\longrightarrow} \lambda_1(L,SR(q_{n,1},Av)^{1/2}).
\end{equation}
Hence $CS(X)\geq  k(n)^{1/4}$.
\end{proof}

Lemma \ref{conway} finishes the proof of Theorem \ref{katz}.
\end{proof}

\section{Symplectic manifolds} \label{symplectic}

We now show that symplectic 4-manifolds with $b_2^+=1$ necessarily have $b_2\leq 10$, as stated in the introduction (note that we always assume symplectic forms to be compatible with the orientation, i.e.~of positive square). 
\begin{lem} \label{symp} Let $X$ be a closed symplectic 4-manifold with $b_2^+=1$ and even intersection form $Q$. Then $Q \cong H$ or $Q \cong H\oplus (-E_8)$. 
\end{lem}
Here $H$ denotes the bilinear form given by $H=\left(\begin{array}{cc} 0 & 1 \\ 1 & 0 \\ \end{array}\right)$ and $E_8$ is a positive-definite, even form of rank 8 associated to the Dynkin diagram of the Lie group $E_8$ (see \cite{GS}).   
\begin{proof} 
If $b_2^-=0$ then $Q =(1)$. If $b_2^->0$ then $Q$ is indefinite and hence of the form $Q=H\oplus (-k)E_8$, since $Q$ is even. It follows that $X$ is minimal because the intersection form does not split off a $(-1)$. If $K^2<0$ then according to a theorem of A.-K.~Liu \cite{Liu}, $X$ is an irrational ruled surface and hence has intersection form $Q=H$ (or $Q=(1)\oplus(-1)$, which is odd). If $K^2=2\chi+3\sigma\geq 0$ then $4b_1+b_2^-\leq 9$ and $b_1=0$ or $b_1=2$, because $1-b_1(X)+b_2^+(X)$ is an even number for every almost complex 4-manifold $X$. If $b_1=0$ then $b_2^-\leq 9$, hence $Q=H$ or $H\oplus (-E_8)$. If $b_1=2$ then $b_2^-\leq 1$, hence $Q=H$. 
\end{proof}

\begin{rem}
It is possible to give a different proof of Proposition \ref{dense} for symplectic manifolds, which relies on a theorem of T.-J.~Li and A.-K.~Liu (\cite{LL3}, Theorem 4). This theorem implies that on a closed 4-manifold $X$ with $b_2^+=1$ which admits a symplectic structure, the set of classes in $H^2(X;\mathbb{R})$ represented by symplectic forms is dense in the positive cone, because it is the complement of at most countably many hyperplanes. If a closed symplectic 4-manifold with $b_2^+=1$ is {\it minimal} (i.e.~there are no symplectic $(-1)$-spheres), then the period map is in fact surjective.
\end{rem}

\section{The $\frac{5}{4}$-conjecture and some examples} \label{5/4} 
The $\frac{5}{4}$-conjecture is a weak analogue of the $\frac{11}{8}$-conjecture which relates the signature and second Betti number of spin 4-manifolds. The main result in this direction is a theorem of M.~Furuta \cite{F}, generalizing work of S.~K.~Donaldson \cite{D1,D2}, that all closed oriented spin 4-manifolds $X$ with $b_2(X)>0$ satisfy the inequality 
\begin{equation}
\frac{5}{4}|\sigma(X)|+2\leq b_2(X),
\end{equation}
where $\sigma(X)$ denotes the signature. C.~Bohr \cite{Bo} then proved a slightly weaker inequality $\frac{5}{4}|\sigma(X)|\leq b_2(X)$ for all 4-manifolds with even intersection form and certain fundamental groups, including all finite and all abelian groups. These are special instances of the following general $\frac{5}{4}$-conjecture. 
\begin{conj} \label{conj} If $X$ is a closed oriented even 4-manifold, then 
\begin{equation} 
\frac{5}{4}|\sigma(X)|\leq b_2(X),
\end{equation}
where $\sigma(X)$ denotes the signature.
\end{conj} 
Here we call a 4-manifold {\it even}, if it has even intersection form.

\begin{rem} J.-H.~Kim \cite{Ki} has proposed a proof of the $\frac{5}{4}$-conjecture. However, some doubts have been raised about the validity of the proof. Therefore, we have chosen to state the result still as a conjecture.
\end{rem} 

\begin{lem} If $X$ is an even 4-manifold with $b_2^+=1$, then the $\frac{5}{4}$-conjecture holds for $X$ if and only if $Q\cong H$ or $Q \cong H\oplus (-E_8)$. 
\end{lem}
\begin{proof} If $X$ is an even 4-manifold with $b_2^+=1$, then $Q\cong H\oplus(-k)E_8$ for some $k\geq 0$. The $\frac{5}{4}$-conjecture is equivalent to $k\leq 1$. 
\end{proof}
In particular, by Lemma \ref{symp}, the $\frac{5}{4}$-conjecture holds for all even symplectic 4-manifolds with $b_2^+=1$. 

There are many examples of 4-manifolds with $b_2^+=1$ where Theorem \ref{coarse} applies, e.g.\ the infinite family of simply-connected pairwise non-diffeomorphic Dolgachev surfaces which are all homeomorphic to $\mathbb{C}P^2\#9\overline{\mathbb{C}P^2}$ (see \cite{GS}). These 4-manifolds are K\"ahler, hence symplectic. There are also recent constructions of infinite families of non-symplectic and pairwise non-diffeomorphic 4-manifolds homeomorphic to $\mathbb{C}P^2\#n\overline{\mathbb{C}P^2}$ for $n\geq 5$ (see \cite{FS2,PSS}). If we take multiple blow-ups of these manifolds, the blow-up formula for the Seiberg-Witten invariants \cite{FS1} shows that the resulting manifolds stay pairwise non-diffeomorphic. Hence we obtain infinite families of symplectic and non-symplectic 4-manifolds $X$ with $n=b_2(X)\rightarrow \infty$, where Theorem \ref{katz} applies.

\subsection*{Acknowledgments.}
I wish to thank D.~Kotschick for raising the question about possible extensions of the work of M.~Katz and many helpful comments.

\end{document}